\documentclass[11pt,reqno]{amsart}

\usepackage{regb-article-preamble}

\usepackage{pgf}
\usepackage{pgffor}
    \usepgfmodule{shapes}
    \usepgfmodule{plot}
\usepackage{tikz}
    \usetikzlibrary{positioning,decorations.pathreplacing,quotes,angles,arrows.meta}

\hypersetup{linktocpage=true,
    colorlinks=true,
    linkcolor=blue,
    citecolor=blue,
    pdftitle={Exceptional projections of sets exhibiting almost dimension conservation},
    pdfauthor={Ryan E. G. Bushling}}

\begin{document}

\begin{abstract}
    We establish a packing dimension estimate on the exceptional sets of orthogonal projections of sets satisfying an \textit{almost dimension conservation} law. In particular, the main result applies to \textit{homogeneous sets} and to certain \textit{graph-directed sets}. Connections are drawn to results of M. Rams and T. Orponen. 
\end{abstract}

\title[]{Exceptional projections of sets exhibiting \\ almost dimension conservation}
\author[]{Ryan E. G. Bushling}
\address{Department of Mathematics \\ University of Washington, Box 354350 \\ Seattle, WA 98195-4350}
\email{reb28@uw.edu}

\maketitle

\section{Introduction} \vs{-0.15}

Let $\Gr(n,k)$ be the Grassmannian of $k$-planes in $\R^n$ equipped with the invariant measure $\gamma_{n,k}$ induced by the action of the orthogonal group. That is, if $\theta_n$ is the Haar measure on $\O(n)$ and $V \in \Gr(n,k)$, then
\begin{equation*}
    \gamma_{n,k}(A) := \theta_n\big( \{ T \in \O(n) \!: T(V) \in A \} \big),
\end{equation*}
where the definition does not depend on the specific choice of $V$. For each $V \in \Gr(n,k)$, we denote by $\pi_V \!: \R^n \to V$ the orthogonal projection of $\R^n$ onto $V$ and write $A_V := \pi_V(A)$ for $A \subseteq \R^n$. The following result of Marstrand, Kaufman, Mattila, and Falconer is essentially the most general bound on the Hausdorff dimension of exceptional sets of orthogonal projections; however, see also \cite{falconer2015sixty} for other results of this sort.

\begin{thm}[Marstrand's Projection Theorem] \label{thm:marstrand-projection}
Let $A \subseteq \R^n$ be analytic.
\begin{enumerate}[label={\normalfont\textbf{(\alph*)}}, itemsep=0pt, topsep=0pt]
    \item If $s \leq \dim A \leq k$, then
    \begin{equation*}
        \dim \, \{ V \in \Gr(n,k) \!: \dim A_V < s \} \leq k(n-k) - (k-s).
    \end{equation*}
    \item If $s \leq k < \dim A$, then
    \begin{equation*}
        \dim \, \{ V \in \Gr(n,k) \!: \dim A_V < s \} \leq k(n-k) - (\dim A - s).
    \end{equation*}
\end{enumerate}
\end{thm}

\subsection{A theorem of Rams}

Theorem 1.1 of Rams \cite{rams2002packing} implies a result in a similar vein that bounds the \textit{packing} dimension of the exceptional set, but only when the set under consideration is highly regular and $k = n-1$. His theorem concerns parametrized families of conformal iterated function systems (IFS), but we state it here more narrowly for parametrized families of similarities. The terminology involved is heavy, so we refer the reader to \S 1 and Definition 4.4 of \cite{rams2002packing} and the \hyperlink{appendix}{Appendix} of this paper for additional background.

\begin{thm}[Rams \cite{rams2002packing}] \label{thm:rams}
Let $V \subset \R^k$ be a bounded open set, and for each $t \in \overline{V}$, let $\big( f_i(\,\cdot\,;t) \big)_{i=1}^N$ be a family of self-similar IFS on $\R^k$ with limit set $K_t$. Assume that each $f_i$ is smooth in all $k$ variables and $k$ parameters, and denote by $\sigma(t)$ the solution to Hutchinson's equation
\begin{equation*}
    \sum_{i=1}^N a_i(t)^{\sigma(t)} = 0,
\end{equation*}
where $a_i(t) \in (0,1)$ is the similarity ratio of $f_i(\,\cdot\,;t)$. If $\big( f_i(\,\cdot\,;t) \big)_{i=1}^N$ satisfies the transversality condition \eqref{eq:transversality}, then
\begin{equation*}
    \dim_P \, \{ u \in \overline{V} \!: \dim K_u \leq s \} \leq s \qquad \forall \, 0 \leq s < \min \left\{ \, k, \, \sup_{t \in \overline{V}} \+ \sigma(t) \right\}.
\end{equation*}
\end{thm}

We adapt our previous notation and let $\rho_e \!: \R^n \to \R^{n-1}$ denote the orthogonal projection onto the hyperplane orthogonal to the vector $e \in \bbs^{n-1}$. Henceforth, ``IFS" without qualification refers to a \textit{self-similar} IFS.

Suppose $K \subset \R^n$ is the limit set of an IFS $(g_i)_{i=1}^N$ that satisfies the strong separation condition (SSC), and for each $e \in \bbs^{n-1}$, let $f_i(\,\cdot\,;e) := \rho_e \circ g_i \circ \rho_e^{-1} \!: \R^{n-1} \to \R^{n-1}$. Then each collection $(f_i(\,\cdot\,;e))_{i=1}^N$ is an IFS in $\R^{n-1}$ with limit set $K_e := \rho_e(K)$. Therefore, $(f_i(\,\cdot\,;e))_{i=1}^N$ is a transverse, smoothly parametrized $(n-1)$-parameter family of IFS on $\R^{n-1}$, so the problem of determining the exceptional set of projections for $K$ is equivalent to determining the exceptional set of the IFS $(f_i(\,\cdot\,;e))_{i=1}^N$.

This setup allows for an application of Rams' theorem to obtain the following.

\begin{prop} \label{prop:rams}
Let $K \subset \R^n$ be the limit set of a self-similar IFS containing no rotations or reflections and satisfying the SSC. Then
\begin{equation*}
    \dim_P \, \{ e \in \bbs^{n-1} \!: \dim K_e \leq s \} \leq s \qquad \forall \, 0 \leq s < \dim K.
\end{equation*}
\end{prop}

Our main result subsumes this as a special case; nevertheless, we include its proof in the \hyperlink{appendix}{Appendix} to shed light on the relationship between his work and our own, and to give a sense of just how strong the transversality condition in Theorem \ref{thm:rams} is.

\subsection{Orponen's packing dimension bound on exceptional sets}

For planar sets, Orponen's result in \cite{orponen2015packing} allows us to forgo any separation conditions in the case that $K$ is self-similar, or to instead assume that $K$ is homogeneous (see \S \ref{ss:homog-and-dc}).

\begin{prop}[Orponen \cite{orponen2015packing}] \label{prop:orponen}
Let $K \subset \R^2$ be homogeneous or self-similar. Then
\begin{equation*}
    \dim_P \, \{ e \in \bbs^1 \!: \dim K_e \leq s \} \leq s \qquad \forall \, 0 \leq s < \dim K.
\end{equation*}
\end{prop}

This result generalizes nicely to higher dimensions in the homogeneous case, but the self-similar case is far more delicate: there is a dichotomy between planar self-similar sets containing dense rotations in the orthogonal group and those that do not, whereas no such dichotomy exists in dimension $n \geq 3$. Peres and Shmerkin \cite{peres2009resonance} showed that the exceptional set is \textit{empty} for self-similar sets in $\R^2$ with dense rotations, leaving only the case of non-dense rotations\textemdash a case that readily reduces to the case of homogeneous sets.

The key to Proposition \ref{prop:orponen} is that linear maps are ``dimension conserving" for homogeneous sets, and the equality of Hausdorff and upper box dimension of homogeneous sets allows one to utilize this principle in a discrete setting. Our main purposes are to generalize Proposition \ref{prop:orponen} to higher dimensions and to weaken the dimension conservation principle to ``\textit{almost} dimension conservation" (cf.\! \S \ref{ss:gd-and-almost-dc}).

\begin{thm} \label{thm:main}
Let $A \subset \R^n$ be a bounded set with $\dim A = \overline{\dim}_B \+ A$. If there exists a set $F \subseteq \Gr(n,k)$ with $t := \dim_P F$ such that $\pi_V$ is almost dimension conserving for $A$ for all $V \in \Gr(n,k) \setminus F$, then
\begin{equation} \label{eq:main}
    \dim_P \, \{ V \in \Gr(n,k) \!: \dim A_V \leq s \} \leq \max\!\big( k(n-k) - (k-s), \+ t \big) \quad \forall \, 0 \leq s < \dim A.
\end{equation}
\end{thm}

For examples of sets exhibiting almost dimension conservation but (to the author's knowledge) not necessarily with dimension conserving projections in the sense of Furstenberg, see \cite{falconer2015dimension} and \cite{farkas2019dimension}.

As a corollary to Theorem \ref{thm:main}, we obtain the following analogue of Proposition \ref{prop:orponen}.

\begin{cor} \label{cor:main}
Let $K \subset \R^n$ be either a homogeneous set or a graph-directed set with one of the following two properties: either (1) its transformation group is finite or (2) the action of the transformation group on $\Gr(n,k)$ has a dense orbit. Then
\begin{equation} \label{eq:cor}
    \dim_P \, \{ V \in \Gr(n,k) \!: \dim K_V \leq s \} \leq k(n-k) - (k-s) \qquad \forall \, 0 \leq s < \dim K.
\end{equation}
\end{cor}

\textit{Proof.} In both cases, $\dim K = \overline{\dim}_B \+ K$. If $K$ is homogeneous, then every linear map (in particular, every orthogonal projection) is (classically) dimension conserving for $K$ by Furstenberg's Theorem \ref{thm:dc}. Likewise, if $K$ is a graph-directed set with finite transformation group, then every linear map is almost dimension conserving by Farkas' Theorem \ref{thm:almost-dc}. Lastly, if $K$ is graph-directed and has a dense orbit in the Grassmannian, then \cite{farkas2019dimension} Theorem 1.15 implies that the exceptional set of projections is empty. \hfill $\square$

If one replaces almost dimension conservation with \textit{bona fide} dimension conservation in Theorem \ref{thm:main}, one can still conclude Corollary \ref{cor:main} using \cite{farkas2019dimension} Theorem 1.3, which implies that graph-directed sets with finite transformation groups are well-approximated from within by homogeneous sets. However, the theorem is stated as is in the interests of generality and future examples of almost dimension conservation.

Our proof of Theorem \ref{thm:main} is similar to Orponen's proof of Proposition \ref{prop:orponen}. Philosophically, dimension conservation affords us partial knowledge of why the dimension of a set may have dropped upon projection onto a subspace: the dimension of the fibers accounts for a ``substantial portion" of the dimension lost\textemdash the dimension did not simply ``vanish"\textemdash and this information enables us to deduce what is happening upstairs in $\R^n$ from what we see downstairs in $\R^k$. While one might hope to prove a still broader restatement of Theorem \ref{thm:main}, \cite{orponen2015packing} Theorem 1.1 shows that \eqref{eq:cor} does not hold in general even for compact sets $K \subset \R^2$. Consequently, any strengthening of Theorem \ref{thm:main} will require \textit{some} mechanism taking the role of dimension conservation.

Kaufman's approach to Theorem \ref{thm:marstrand-projection}(a) relies on potential theory, and our proof of Theorem \ref{thm:main} can be understood as a discretization thereof. Theorem \ref{thm:marstrand-projection}(b) begs the following improvement to Theorem \ref{thm:main}; however, as Theorem \ref{thm:marstrand-projection}(b) is not known to be provable without use of the Fourier transform, it is probable that any modification of our proof leading to such an extension will also need to emulate the role played by the Fourier transform.

\textbf{Conjecture.} \textit{Under the hypotheses of Theorem \ref{thm:main}, if $\dim A > k$, then
\begin{equation*}
    \dim_P \, \{ V \in \Gr(n,k) \!: \dim A_V < s \} \leq \max\!\big( k(n-k) - (\dim A - s), \+ t \big) \quad \forall \, k < s < \dim A.
\end{equation*}}

\subsection{Outline}

In \S \ref{s:definitions}, we lay down some preliminaries on Hausdorff and packing dimensions, followed by the requisite background on dimension conservation, almost dimension conservation, homogeneous sets, and graph-directed sets. A lemma on the sizes of discrete sets of points in $\Gr(n,k)$ follows in \S \ref{s:gr(n,k)}, and this is used in \S \ref{s:main-theorem} in the proof of Theorem \ref{thm:main}. Finally, the proof of Proposition \ref{prop:rams} applying Rams' Theorem \ref{thm:rams} to orthogonal projections of self-similar sets appears in the \hyperlink{appendix}{Appendix}.

\section{Definitions and notation} \label{s:definitions}

Throughout this paper, the relation $a \lesssim b$ denotes non-strict inequality up to a positive multiplicative constant and the relation $a \sim b$ indicates that both $a \lesssim b$ and $b \lesssim a$. For clarity, the parameters on which the implicit constant depends will sometimes be written as subscripts or stated explicitly.

\subsection{Measures and dimensions}

Since we will be working in both Euclidean space and the Grassmannian $\Gr(n,k)$\textemdash a metric space with metric $d(V,W) := \| \pi_V - \pi_W \|$\textemdash the following discussion is framed in a metric space context. However, little is lost by taking the metric space to be an open subset of $\R^n$.

Let $(X,d)$ be a separable metric space and $2^X$ its power set, and let $|F|$ be the diameter of a set $F \in 2^X$. The Carath\'eodory construction yields a family of functions $\mathcal{H}_\delta^s \!: 2^X \to [0,\infty]$, $\delta \in [0,\infty]$, defined by
\begin{equation*}
    \mathcal{H}_\delta^s(A) := \inf \left\{ \sum_{i=1}^\infty |F_i|^s \!: F_i \in 2^X, \ |F_i| \leq \delta, \ A \subseteq \bigcup_{i=1}^\infty F_i \right\}.
\end{equation*}
The function $\mathcal{H}_\infty^s$, called the \textit{$s$-dimensional Hausdorff content} on $X$, will be of particular importance below. The resulting Carath\'eodory measure
\begin{equation*}
    \mathcal{H}^s(A) := \lim_{\delta \to 0} \mathcal{H}_\delta^s(A) = \sup_{\delta > 0} \mathcal{H}_\delta^s(A)
\end{equation*}
is called the \textit{$s$-dimensional Hausdorff measure} on $X$. For each $A \subseteq X$, there is a unique $s \in [0,\infty]$ (in fact, $s \in [0,n]$ when $X$ is a smooth $n$-manifold) with the following property: for all $r < s < t$, we have
\begin{equation*}
    0 = \mathcal{H}^t(A) \leq \mathcal{H}^s(A) \leq \mathcal{H}^r(A) = \infty,
\end{equation*}
or, equivalently,
\begin{equation*}
    s = \sup \, \big\{ t \in [0,\infty) \!: \mathcal{H}^t(A) > 0 \big\} = \inf \, \big\{ r \in [0,\infty) \!: \mathcal{H}^r(A) = 0 \big\}.
\end{equation*}
We write $\dim A := s$ and call this number the \textit{Hausdorff dimension} of $A$.

Perhaps the most useful nontrivial property of Hausdorff dimension that we use is the following, as it entails that Hausdorff content is sufficient to determine Hausdorff dimension.

\begin{prop} \label{prop:hdorff-content}
Let $A \subseteq X$. Then $\mathcal{H}^s(A) > 0$ if and only if $\mathcal{H}_\infty^s(A) > 0$.
\end{prop}

Hausdorff measure\textemdash hence, Hausdorff dimension\textemdash is defined in terms of covers by \textit{arbitrary} sets of diameter \textit{at most} $\delta$. However, one can still recover information about Hausdorff dimension using a smaller family of covers, namely, covers by \textit{balls} with radius \textit{equal to} $\delta$. For each bounded set $A \subseteq X$ and each $\delta > 0$, let
\begin{equation} \label{eq:covering-number}
    N(A,\delta) := \min \left\{ k \in \Z_+ \!: \exists \, x_i \in X \text{ s.t. } A \subseteq \bigcup_{i=1}^k B(x_i,\delta) \right\},
\end{equation}
where (for definiteness) we take the balls to be closed. We define the \textit{upper box dimension} of a set in $X$, also called its \textit{upper Minkowski dimension}, by
\begin{equation} \label{eq:upper-box-dim}
    \begin{aligned}
    \overline{\dim}_B \+ A :\!&= \sup \left\{ t \in [0,\infty) \!: \limsup_{\delta \downarrow 0} N(A,\delta) \+ \delta^t > 0 \right\} \\
    &= \inf \left\{ r \in [0,\infty) \!: \limsup_{\delta \downarrow 0} N(A,\delta) \+ \delta^r = 0 \right\}.
    \end{aligned}
\end{equation}
Comparing the admissible covers of $A$ in the definitions of Hausdorff and upper box dimensions yields the inequality
\begin{equation*}
    \dim A \leq \overline{\dim}_B \, A.
\end{equation*}
We will frequently deal with sets $A$ for which equality holds.

We now turn to packing dimension. It is perhaps more natural to define this in terms of packing measure, but, to streamline the exposition, we present an alternative characterization\textemdash the only one we will use in the sequel. The \textit{packing dimension} of any subset $A \subseteq X$ is given by
\begin{equation*}
    \dim_P A := \inf \left\{ \sup_{i \in \Z_+} \overline{\dim}_B \+ A_i \!: A = \bigcup_{i=1}^\infty A_i, \, |A_i| < \infty \right\}.
\end{equation*}
An important feature of packing dimension that upper box dimension lacks is \textit{countable stability}:
\begin{equation*}
    \dim_P \bigcup_{i=1}^\infty A_i = \sup_{i \in \Z_+} \dim_P A_i. \vs{0.15}
\end{equation*}

\subsection{Dimension conservation and homogeneous sets} \label{ss:homog-and-dc}

As a prelude to the introduction of ``almost dimension conservation," we discuss the relevant terminology from Furstenberg's paper \cite{furstenberg2008ergodic}.

A Lipschitz function $f \!: \R^n \to \R^m$ is said to be \textit{dimension conserving} (or \textit{DC}) for a set $A \subseteq \R^n$ if there exists $\Delta \geq 0$ such that
\begin{equation*}
    \Delta + \dim \big\{ y \in f(A) \!: \dim(A \cap f^{-1}(y)) \geq \Delta \big\} \geq \dim A,
\end{equation*}
where $\dim \varnothing := -\infty$. Heuristically, $f$ is DC for $A$ if the dimension of its fibers over $A$ is ``complementary" to the dimension of $f(A)$: it is a sort of ``rank-nullity condition." The pathological Example 7.8 of \cite{falconer2014fractal} shows that even the projection of a product set onto the coordinate axes may radically fail to be DC for that set.

The \textit{Hausdorff metric} on the class $\mathcal{K}$ of nonempty compacta in $\R^n$ is defined by
\begin{equation*}
    \rho(H,K) := \inf \, \{ \eps \geq 0 \!: H \subseteq K_\eps \text{ and } K \subseteq H_\eps \},
\end{equation*}
$H,K \in \mathcal{K}$, where $A_\eps$ is the closed $\eps$-neighborhood of $A$. With the Hausdorff metric, $\mathcal{K}$ is a complete metric space.

We now define the archetypal class of compacta in the study of dimension conservation. Scaling and translating a set does not affect the dimension of the projection of a set in any direction, so we assume without loss of generality that $K \subseteq [0,1]^n$. A closed set $K' \subset [0,1]^n$ is called a \textit{miniset} of $K$ if there exists an expanding homothety $\varphi(x) = rx + b$ ($|r| \geq 1$) such that $K' \subseteq \varphi(K)$. A closed set $K'' \subseteq [0,1]^n$ is called a \textit{microset} of $K$ if there exists a sequence $(K_i')_{i=1}^\infty$ of minisets of $K$ converging to $K''$ in the Hausdorff metric: $\rho(K_i',K'') \to 0$. Finally, $K$ is said to be \textit{homogeneous} if all its microsets are minisets; that is, the class of minisets of $K$ is a closed in $\mathcal{K}$.

Loosely, $K$ is homogeneous if it looks the same at all scales: even if the minisets $K_i'$ must be contained in larger and larger expansions of $K$ as $i \to \infty$ (meaning they resemble smaller and smaller subsets of $K$), there still exists a scale on which the limiting set $K''$ coincides with a subset of $K$ at that scale.

Appreciation for the definition of ``dimension conserving" will be important for understanding the sequel. On the other hand, the technical definition of a homogeneous set is not strictly necessary, as the only two properties we require are the following (cf.\! \cite{furstenberg2008ergodic} p.\! 407 and Theorem 6.2).

\begin{prop} \label{prop:homogeneous-dimension}
If $K$ is homogeneous, then $\dim K = \overline{\dim}_B \+ K$.
\end{prop}

\begin{thm}[Furstenberg \cite{furstenberg2008ergodic}] \label{thm:dc}
If $K \subset \R^n$ is homogeneous and $f \!: \R^n \to \R^m$ is linear, then $f$ is DC for $K$.
\end{thm}

In particular, Theorem \ref{thm:dc} implies that every projection map is DC for $K$.

\subsection{Almost dimension conservation and graph-directed sets} \label{ss:gd-and-almost-dc}

There is a natural weakening of the notion of dimension conservation due to \cite{farkas2019dimension}. We call a Lipschitz function $f \!: \R^n \to \R^m$ \textit{almost dimension conserving} for a set $A \subseteq \R^n$ if there exists $\Delta \geq 0$ such that, for every $\eps > 0$,
\begin{equation} \label{eq:almost-DC}
    \Delta + \dim \big\{ y \in f(A) \!: \dim(A \cap f^{-1}(y)) \geq \Delta - \eps \big\} \geq \dim A,
\end{equation}
where $\dim \varnothing := -\infty$ as before. See also \cite{falconer2015dimension} for the more restrictive concept of \textit{weak} dimension conservation. In \S \ref{s:main-theorem}, $A$ is fixed and we let $\Delta'(V)$ denote the set of all $\Delta \geq 0$ such that \eqref{eq:almost-DC} holds with $f = \pi_V$.

The primary motivation for working with almost dimension conservation is its applicability to a wider class of sets, the most notable of which we describe as follows. Let $(\mathcal{E}, \mathcal{V})$ be a directed graph with vertices $1, ..., N$, where edges starting and ending at the same point are allowed. For all $i,j \in \mathcal{V}$, let $\mathcal{E}_{i,j}$ be the set of all edges from $i$ to $j$ and $\mathcal{E}_{i,j}^\ell$ the set of all paths of length $\ell$ from $i$ to $j$. We assume that $(\mathcal{E}, \mathcal{V})$ is \textit{transitive} (or \textit{strongly connected}), i.e., that for all $i,j \in \mathcal{V}$, at least one of the $\mathcal{E}_{i,j}^\ell$ is nonempty. Given a family $(g_e)_{e \in \mathcal{E}}$ of contracting similarities, there exists a unique tuple $(K_1, ..., K_N)$ of compacta with the invariance property
\begin{equation*}
    K_i = \bigcup_{j=1}^N \bigcup_{e \in \mathcal{E}_{i,j}} g_e(K_j)
\end{equation*}
for each $i$. The system $(g_e)_{e \in \mathcal{E}}$ is called a \textit{graph-directed IFS}, the tuple $(K_1, ..., K_N)$ is called its \textit{attractor}, and each $K_i$ is called a \textit{graph-directed set}. Finally, if we write $g_e(x) = a_e \+ T_e(x) + b_e$, where $a_e \in (0,1)$ are the similarity ratios, $T_e \in \O(n)$ are orthogonal transformations, and $b_e \in \R^n$ are translation vectors, then the group $\mathcal{T}_i \leq \O(n)$ generated by the set $\{ T_{e_1} \circ \cdots \circ T_{e_\ell} \in \O(n) \!: (e_1, ..., e_\ell) \in \mathcal{E}_{i,i}^\ell, \, \ell \in \Z_+ \}$ is called the \textit{transformation group} of $K_i$. Of particular importance are the case that $\mathcal{T}_i$ is finite and the case that some (hence, every) orbit in $\Gr(n,k)$ under the action of $\mathcal{T}_i$ is dense. What happens between these two extremes, while interesting, is not well-understood (see \cite{shmerkin2015projections} \S 7.2).

Analogous to Proposition \ref{prop:homogeneous-dimension} and Theorem \ref{thm:dc} are the following. They are the only properties of graph-directed sets important for our purposes.

\begin{prop} \label{prop:graph-dimension}
If $K$ is a graph-directed set, then $\dim K = \overline{\dim}_B \+ K$.
\end{prop}

\begin{thm}[Farkas \cite{farkas2019dimension}] \label{thm:almost-dc}
Let $(K_1, ..., K_\ell)$ be the attractor of a graph-directed IFS in $\R^n$. If $i \in \{ 1, ..., \ell \}$ is any vertex such that the transformation group of $K_i$ is finite and if $f \!: \R^n \to \R^m$ is linear, then $f$ is almost DC for $K_i$.
\end{thm}

Thus, under the hypotheses of Theorem \ref{thm:almost-dc}, every projection map is almost DC for $K_i$.

\section[Counting points on Gr(n,k)]{Counting points on \texorpdfstring{$\Gr(n,k)$}{Gr(n,k)}} \label{s:gr(n,k)}

We will have occasion to ask, given two points $x,y \in \R^n$ and a $\delta$-separated set $E \subset \Gr(n,k)$, for how many $V \in E$ we have $\| \pi_V(x) - \pi_V(y) \| \leq c \+ \delta$, where $c > 0$ is given. The following lemma\textemdash a discrete version of a key step in the proof of Marstrand's projection theorem\textemdash addresses this question in greater generality.

\begin{lem} \label{lem:discrete-proj}
Let $x \in \R^n \setminus \{ 0 \}$ and $\delta_1, \delta_2 \in (0,1]$ and let $E \subset \Gr(n,k)$ be $\delta_2$-separated. Then
\begin{equation*}
    \card \big\{ V \in E \!: \| \pi_V(x) \| \leq \delta_1 \big\} \+\lesssim_{n,k}\+ \delta_1^k \+ \delta_2^{-k(n-k)} \+ \| x \|^{-k}.
\end{equation*}
\end{lem}

\textit{Proof.} Since the invariant measure $\gamma_{n,k}$ and the $k(n-k)$-dimensional Hausdorff measure $\mathcal{H}^{k(n-k)}$ are both uniformly distributed measures on $\Gr(n,k)$, they are equal up to a constant. Consequently, $r^{k(n-k)} \lesssim \gamma_{n,k}(B(V,r))$ for all $V \in \Gr(n,k)$ and $r \in (0,1]$, and it follows from the separation hypothesis on $E$ that
\begin{equation*}
    \delta_2^{k(n-k)} \card \big\{ V \in E \!: \| \pi_V(x) \| \leq \delta_1 \big\} \+\lesssim\+ \gamma_{n,k}\big( \big\{ V \in \Gr(n,k) \!: \| \pi_V(x) \| \leq \delta_1 \big\} \big).
\end{equation*}
Lemma 3.11 of \cite{mattila1995geometry} states that
\begin{equation*}
    \gamma_{n,k}\big( \big\{ V \in \Gr(n,k) \!: \| \pi_V(x) \| \leq \delta_1 \big\} \big) \lesssim_{n,k} \delta_1^k \+ \| x \|^{-k},
\end{equation*}
and combining this inequality with the previous gives
\begin{equation*}
    \delta_2^{k(n-k)} \card \big\{ V \in E \!: \| \pi_V(x) \| \leq \delta_1 \big\} \+\lesssim\+ \delta_1^k \+ \| x \|^{-k}.
\end{equation*}
Dividing through by $\delta_2^{k(n-k)}$ yields the desired inequality. \hfill $\square$

The question at the start of this section is answered by replacing $x$ with $x-y$ and applying the linearity of $\pi_V$. Curiously, the proof of Theorem \ref{thm:main} (and, in particular, of Lemma \ref{lem:second-DC-lem}) does not seem to benefit from the full generality of Lemma \ref{lem:discrete-proj}: discretizing in the Grassmannian and in Euclidean space at the same scale (i.e., letting $\delta_1 \sim \delta_2$) does not seem to affect the strength of the conclusion.

\section[Proof of Theorem 1.5]{Proof of Theorem \ref{thm:main}} \label{s:main-theorem}

Theorem \ref{thm:main} will follow readily from the following two lemmas. The first is essentially true by definition, but it quantitatively formalizes the idea that, if $K$ is homogeneous and the dimension of $K_V$ is small, then the dimensions of the fibers of $\pi_V$ over $K$ must be large.

\begin{lem} \label{lem:first-DC-lem}
Let $A \subseteq \R^n$. If there exists a set $F \subset \Gr(n,k)$ with $t := \dim_P F$ such that $\pi_V$ is almost DC for $A$ for all $V \in \Gr(n,k) \setminus F$, then
\begin{align*}
    & \dim_P \big\{ V \in \Gr(n,k) \!: \dim A_V \leq s \big\} \\
    &\hs{1.2} \leq \max\!\Big( t, \, \dim_P \big\{ V \in \Gr(n,k) \setminus F \!: \Delta + s \geq \dim A \ \forall \Delta \in \Delta'(V) \big\} \Big).
\end{align*}
for all $s \geq 0$.
\end{lem}
\textit{Proof.} Let $V \in \Gr(n,k) \setminus F$ be such that $\dim A_V \leq s$ and let $\Delta \in \Delta'(V)$. Then, by the definition of almost dimension conservation,
\begin{equation*}
    \Delta + \dim \big\{ y \in A_V \!: \dim\!\big( A \cap \pi_V^{-1}(y) \big) \geq \Delta - \eps \big\} \geq \dim A
\end{equation*}
for every $\eps > 0$. It follows from the monotonicity of dimension that
\begin{equation*}
    \Delta + s \geq \Delta + \dim A_V \geq \dim A.
\end{equation*}
This proves the inclusion
\begin{equation*}
    \big\{ V \in \Gr(n,k) \setminus F \!: \dim A_V \leq s \big\} \subseteq \big\{ V \in \Gr(n,k) \setminus F \!: \Delta + s \geq \dim A \ \forall \Delta \in \Delta'(V) \big\},
\end{equation*}
whence
\begin{align*}
    & \dim_P \big\{ V \in \Gr(n,k) \!: \dim A_V \leq s \big\} \\
    &\hs{1} \leq \max\!\Big( \dim_P F, \, \dim_P \big\{ V \in \Gr(n,k) \setminus F \!: \dim A_V \leq s \big\} \Big) \\
    &\hs{1} \leq \max\!\Big( t, \, \dim_P \big\{ V \in \Gr(n,k) \setminus F \!: \Delta + s \geq \dim A \ \forall \Delta \in \Delta'(V) \big\} \Big). \qedtag
\end{align*}

\begin{lem} \label{lem:second-DC-lem}
Let $A \subset \R^n$ be a bounded set with $\gamma := \dim A = \overline{\dim}_B \+ A$, and denote
\begin{equation*}
    E_s := \big\{ V \in \Gr(n,k) \!: \exists \+ \Delta \in \Delta'(V) \text{ s.t. } \Delta + s \geq \gamma \big\}.
\end{equation*}
Then
\begin{equation} \label{eq:second-DC-lem}
    \dim_P E_s \leq k(n-k) - (k-s) \qquad \forall \+ 0 \leq s < \gamma.
\end{equation}
\end{lem}

As the proof of this lemma is much more involved, we first show how it (almost immediately) implies Theorem \ref{thm:main}.

\define{Proof of Theorem \ref{thm:main}.} Let $A \subset \R^n$ be a bounded set with $\gamma := \dim A = \overline{\dim}_B \+ A$ whose orthogonal projections $\pi_V$ are almost DC for $V \in \Gr(n,k)$ outside a set $F$ of packing dimension $t$. Then $A$ satisfies the hypotheses of Lemmas \ref{lem:first-DC-lem} and \ref{lem:second-DC-lem}, whence
\begin{align*}
    & \dim_P \big\{ V \in \Gr(n,k) \!: \dim A_V \leq s \big\} \\
    &\hs{0.4} \leq \max\!\Big( t, \, \dim_P \big\{ V \in \Gr(n,k) \setminus F \!: \Delta \geq \gamma - s \ \forall \Delta \in \Delta'(V) \big\} \Big) \\
    &\hs{0.4} \leq \max\!\Big( t, \, \dim_P \big\{ V \in \Gr(n,k) \!: \exists \+ \Delta \in \Delta'(V) \text{ s.t. } \Delta \geq \gamma - s \big\} \Big) \\
    &\hs{0.4} \leq \max\!\big( t, \+ k(n-k) - (k-s) \big). \qedtag
\end{align*}

We conclude by proving Lemma \ref{lem:second-DC-lem}, wherein the dimension conservation hypothesis and the geometry of the Grassmannian come into play.

\define{Proof of Lemma \ref{lem:second-DC-lem}.} The $s = 0$ case follows from the $s > 0$ case by letting $s \downarrow 0$, so we assume that $s > 0$.

\textsc{Step 1.} We begin by reducing the lemma to the following claim: that
\begin{equation} \label{eq:main-reduction}
    \overline{\dim}_B \, E_{s,\eps} \leq k(n-k) - (k-s) + 3\eps
\end{equation}
for all sufficiently small $\eps > 0$, where
\begin{equation*}
    E_{s,\eps} := \Big\{ V \in \Gr(n,k) \!: \exists \Delta \geq \gamma - s \text{\+ s.t.~} \mathcal{H}_\infty^{\gamma - \Delta - \eps}\big( \big\{ y \in \R^k \!: \mathcal{H}_\infty^{\Delta - \eps} \big( A \cap \pi_V^{-1}(y) \big) > \eps \big\} \big) > \eps \Big\}.
\end{equation*}
To prove the validity of this reduction, suppose that $\Delta \geq \gamma - s$ for some $\Delta \in \Delta'(V)$. (In particular, $\pi_V$ is almost DC for $A$.) Then, for every $\eps > 0$,
\begin{equation*}
    \dim \big\{ y \in \R^k \!: \dim\!\big( A \cap \pi_V^{-1}(y) \big) \geq \Delta - \tfrac{\eps}{2} \big\} > \gamma - \Delta - \eps.
\end{equation*}
Consequently, if $0 < \eps < \gamma - s \leq \Delta$, the set on the left-hand side of this inequality has infinite $(\gamma - \Delta - \eps)$-dimensional Hausdorff measure and, in particular, positive $(\gamma - \Delta - \eps)$-dimensional Hausdorff content:
\begin{equation*}
    \mathcal{H}_\infty^{\gamma - \Delta - \eps}\big( \big\{ y \in \R^k \!: \dim\!\big( A \cap \pi_V^{-1}(y) \big) \geq \Delta - \tfrac{\eps}{2} \big\} \big) > 0.
\end{equation*}
Similarly, $\dim\!\big( A \cap \pi_V^{-1}(y) \big) \geq \Delta - \tfrac{\eps}{2}$ implies $\mathcal{H}_\infty^{\Delta - \eps} \big( A \cap \pi_V^{-1}(y) \big) > 0$, so
\begin{equation*}
    \mathcal{H}_\infty^{\gamma - \Delta - \eps}\big( \big\{ y \in \R^k \!: \mathcal{H}_\infty^{\Delta - \eps} \big( A \cap \pi_V^{-1}(y) \big) > 0 \big\} \big) > 0.
\end{equation*}
Writing
\begin{equation*}
    \big\{ y \in \R^k \!: \mathcal{H}_\infty^{\Delta - \eps} \big( A \cap \pi_V^{-1}(y) \big) > 0 \big\} = \bigcup_{m=1}^\infty \big\{ y \in \R^k \!: \mathcal{H}_\infty^{\Delta - \eps} \big( A \cap \pi_V^{-1}(y) \big) > m^{-1} \big\},
\end{equation*}
we can by Proposition \ref{prop:hdorff-content} find $m \in \Z_+$ such that
\begin{equation*}
    \mathcal{H}_\infty^{\gamma - \Delta - \eps}\big( \big\{ y \in \R^k \!: \mathcal{H}_\infty^{\Delta - \eps} \big( A \cap \pi_V^{-1}(y) \big) > m^{-1} \big\} \big) > 0.
\end{equation*}
The left-hand side is decreasing in $\eps$ and increasing in $m$, so we may adjust these parameters accordingly so that $m^{-1} = \eps$ for simplicity. Indeed, for $\eps > 0$ sufficiently small, we will have
\begin{equation*}
    \mathcal{H}_\infty^{\gamma - \Delta - \eps}\big( \big\{ y \in \R^k \!: \mathcal{H}_\infty^{\Delta - \eps}\big( A \cap \pi_V^{-1}(y) \big) > \eps \big\} \big) > \eps,
\end{equation*}
from which we conclude that $V \in E_{s,\eps}$. In particular, we can write
\begin{equation*}
    E_s \subseteq \bigcup_{m=N}^\infty E_{s,m^{-1}}
\end{equation*}
for any $N \in \Z_+$. It then follows from our definition of packing dimension that \eqref{eq:main-reduction} implies \eqref{eq:second-DC-lem}, so we set out to prove \eqref{eq:main-reduction}.

\textsc{Step 2.} Let $0 < \eps < \gamma - s$. We discretize the problem and define a family of relations, indexed by $V \in \Gr(n,k)$, that relate distant points $x,y \in A$ whose projections $\pi_V(x), \pi_V(y)$ are close.

Let $\gamma' > \gamma$, $d := (\gamma - s - \eps)^{-1}$, and
\begin{equation*}
    \delta < \eta := \frac{\eps^d}{n^{1/2} \+ 2^{d+2} \+ (2^{n-k} + 1)^d}.
\end{equation*}
(The significance of this requirement on $\delta$ will become apparent later.) By the definition of upper box dimension, there exists a finite subset $A' \subseteq A$ such that $\card A' \lesssim \delta^{-\gamma'}$ and
\begin{equation*}
    A \subseteq \bigcup_{x \in A'} B(x,\delta).
\end{equation*}
For each $V \in \Gr(n,k)$, let $\mathcal{T}_V$ be the family of \textit{$\delta$-fat $(n-k)$-planes} of the following form:
\begin{equation*}
    \pi_V^{-1}\!\left( \prod_{i=1}^k \big[ j_i \delta, (j_i + 1) \delta \big) \right) \subset \R^n,
\end{equation*}
$j_1, ..., j_k \in \Z$. (For succinctness, we will call these ``fat planes" or ``elements of $\mathcal{T}_V$.") These are half-open neighborhoods of fibers of $\pi_V$ over points of the lattice $\big( \Z + \tfrac{1}{2} \big)^k \delta$, and their disjoint union is all of $\R^n$. We define relations $\backsim_{\scriptscriptstyle{V}}$ on $\R^n$ by
\begin{equation} \label{eq:sim-e}
    \begin{aligned}
    x \backsim_{\scriptscriptstyle{V}} y \quad \iff \quad & \| x-y \| > 2\eta = \frac{\eps^d}{2^{d+1} (2^{n-k} + 1)^d} \quad \text{and} \quad \exists \, T \in \mathcal{T}_V \quad \text{s.t.} \\
    & B(x,\delta) \cap T \neq \varnothing \quad \text{and} \quad B(y,\delta) \cap T \neq \varnothing.
    \end{aligned}
\end{equation}
This states that $\backsim_{\scriptscriptstyle{V}}$ relates points of $\R^n$ that are not too close to each other, but that nevertheless belong to the same fat plane, adjacent fat planes, or fat planes with a common neighboring fat plane. In particular, although the points are fairly distant from each other, their projections onto $V$ are quite close.

\textsc{Step 3.} Let $E' \subseteq E_{s,\eps}$ be any $\delta$-separated subset and define the energy of $E'$ by
\begin{equation} \label{eq:energy}
    \mathcal{E} := \sum_{V \in E'} \card \big\{ (x,y) \in (A')^2 \!: x \backsim_{\scriptscriptstyle{V}} y \big\}.
\end{equation}
We use this energy to bound $\card E'$ and, in turn, $\overline{\dim}_B \, E_{s,\eps}$.

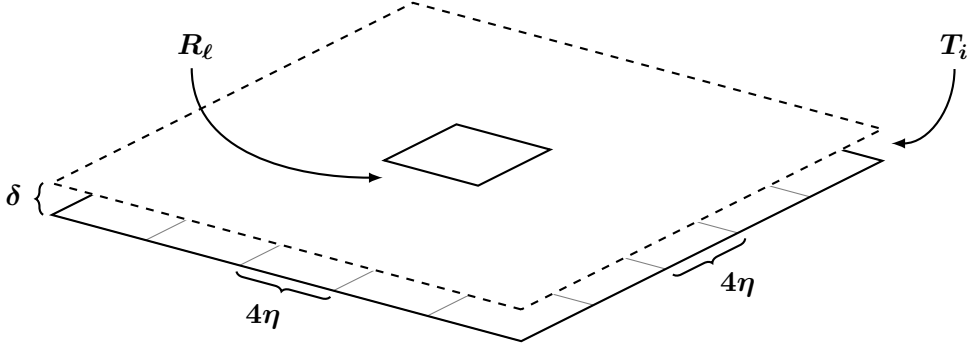
\begin{figure}[ht]

\vs{0.3}

\centering

\begin{tikzpicture}[every node/.style={minimum size=1cm},on grid,scale=1.2]
		

\begin{scope}[yshift=-120,every node/.append style={yslant=0.5,xslant=-1.3},
        yslant=0.5,xslant=-1.3]

    \draw[step=8mm,thin,gray] (0,0) grid (4,4); 
    \draw[black,thick] (1.6,1.6) rectangle (2.4,2.4);
    \draw[black,thick] (0,0) rectangle (4,4); 

    \draw[decorate, decoration={brace, amplitude=1mm}, thick] (-0.1,1.55) -- (-0.1,2.35);
    \draw[decorate, decoration={brace, mirror, amplitude=1mm}, thick] (1.55,-0.1) -- (2.35,-0.1);

    \coordinate (p) at (-0.65,1.7);
    \coordinate (q) at (1.65,-0.575);

    \coordinate (a) at (-0.03,4.05);

\end{scope} 

\draw[thick,white] (-1.535,-2.58) -- (-0.72,-2.18);
\draw[thick,dashed] (-1.535,-2.58) -- (-0.72,-2.18);
\draw[thick,white] (-0.72,-2.18) -- (0.335,-2.46);
\draw[thick,dashed] (-0.72,-2.18) -- (0.335,-2.46);

\draw[thick] (-1.535,-2.58) -- (-1.535,-2.23);
\draw[thick] (-0.475,-2.85) -- (-0.475,-2.50);
\draw[thick] (0.336,-2.45) -- (0.336,-2.104);
\draw[thick,dashed] (-0.72,-2.18) -- (-0.72,-1.83);

\begin{scope}[yshift=-110,every node/.append style={yslant=0.5,xslant=-1.3},
        yslant=0.5,xslant=-1.3]

    \fill[white,fill opacity=0.7] (0,0) rectangle (4,4);
    \draw[white,thick] (0,0) rectangle (4,4); 
    \draw[black,thick,dashed] (0,0) rectangle (4,4); 

    \coordinate (b) at (-0.03,4.05);

    \coordinate (Rla) at (2.2,4.5);
    \coordinate (Rlb) at (2.53,4.72);

    \pgfkeys{/pgf/number format/.cd, fixed, zerofill, precision =1}

\end{scope}

\draw[black,thick,yshift=-100,every node/.append style={yslant=0.5,xslant=-1.3},
        yslant=0.5,xslant=-1.3] (1.16,1.26) rectangle (1.96,2.06);

\draw[-latex,thick] (4.8,-1.2) to[out=-90,in=0] (4.1,-2.025);
\node at (4.8,-0.95) {$\bm{T_i}$};

\draw[decorate, decoration={brace, amplitude=1mm}, thick] (a) -- (b);

\draw[-latex,thick] (Rla) to[out=270,in=180] (-1.535,-2.405);
\node at (Rlb) {$\bm{R_\ell}$};

\node at (p) {$\bm{4\eta}$};
\node at (q) {$\bm{4\eta}$};

\node at (-5.625,-2.6) {$\bm{\delta}$};

\end{tikzpicture}

    \caption{A ``checkerboard square" $R_\ell$ in the fat plane $T_i$ in the case $(n,k) = (3,1)$.}
    \label{fig:checkerboard}
    
    \vs{0.15}
    
\end{figure}

To obtain an upper bound, note that, given $x,y \in A'$, the number of $k$-planes $V \in E'$ such that $x \backsim_{\scriptscriptstyle{V}} y$ is $\lesssim \delta^{-k(n-k)+k} \| x-y \|^{-k}$ by Lemma \ref{lem:discrete-proj}. Hence, for a fixed $x \in A'$, the bound $\card A' \lesssim_{n,\gamma} \delta^{-\gamma'}$ and the trivial estimate $\eta^{-k} \sim_{\eps,n,k} 1$ give
\begin{align*}
    \sum_{V \in E'} & \card \big\{ y \in A' \!: x \backsim_{\scriptscriptstyle{V}} y \big\} = \sum_{y \in A'} \card \big\{ V \in E' \!: x \backsim_{\scriptscriptstyle{V}} y \big\} \\
    & \lesssim \delta^{-k(n-k)+k} \left( \max_{\substack{y \in A', V \in E': \\ x \backsim_{\scriptscriptstyle{V}} y}} \| x-y \|^{-k} \right) \card A' \\[0.15cm]
    &\leq \delta^{-k(n-k)+k} (2\eta)^{-k} \delta^{-\gamma'} \sim \delta^{-k(n-k)+k-\gamma'},
\end{align*}
where the implicit constants depend on $n$, $k$, $\gamma$, $s$, and $\eps$, but not $\gamma'$ or $\delta$. Summing over all $x \in A'$ gives
\begin{equation*}
    \mathcal{E} \lesssim \delta^{-k(n-k)+k-\gamma'} \card A' \lesssim \delta^{-k(n-k)+k-2\gamma'}.
\end{equation*}
To place a lower bound on $\mathcal{E}$, we estimate the individual terms in the sum \eqref{eq:energy}. Let $V \in E' \subseteq E_{s,\eps}$ and $\Delta \in \Delta'(V)$. Unwinding the definition of $E_{s,\eps}$, we see that there exist $j \gtrsim \delta^{\Delta - \gamma + \eps}$ fat planes $T_i \in \mathcal{T}_V$ and points $y_i \in \pi_V(T_i)$, $i = 1, ..., j$, with the following property: if $W_i := \pi_V^{-1}(y_i)$ denotes the $(n-k)$-plane contained in $T_i$ that ``passes through" $y_i$, then
\begin{equation} \label{eq:plane-Hausdorff-content}
    \mathcal{H}_\infty^{\Delta - \eps}(A \cap W_i) > \eps.
\end{equation}
To count the number of relations $x \backsim_{\scriptscriptstyle{V}} y$ that hold on $A'$, we checkerboard each fat plane $T_i \in \mathcal{T}_V$ with boxes or ``checkerboard squares"
\begin{equation*}
    R = T_i \cap \pi_{V^\perp}^{-1}\!\left( \prod_{\ell=1}^{n-k} \big[ 4 i_\ell \eta, 4 (i_\ell + 1) \eta \big) \right),
\end{equation*}
$i_1, ..., i_{n-k} \in \Z$ (see Figure \ref{fig:checkerboard}). Recalling that we chose $\delta < \eta$, we see that
\begin{equation*}
    |R| = \left( \sum_{\ell=1}^{n-k} (4\eta)^2 + \sum_{\ell=n-k+1}^n \delta^2 \right)^{\!1/2} < \left( \sum_{\ell=1}^n (4\eta)^2 \right)^{\!1/2} = n^{1/2} \cdot 4\eta = \frac{\eps^d}{2^d \+ (2^{n-k} + 1)^d}
\end{equation*}
and, consequently, that
\begin{equation} \label{eq:square-Hausdorff-content}
    \mathcal{H}_\infty^{\Delta - \eps}(R) \leq |R|^{\Delta - \eps} < \frac{\eps}{2(2^{n-k} + 1)},
\end{equation}
per our choice of $d$. It then follows from \eqref{eq:plane-Hausdorff-content} and \eqref{eq:square-Hausdorff-content} that, for any choice of squares $R_1, ...,$ $R_{2^{n-k} + 1}$,
\begin{equation*}
    \mathcal{H}_\infty^{\Delta - \eps}\!\left( (A \cap W_i) \setminus \bigcup_{\ell=1}^{2^{n-k} + 1} R_\ell \right) > \frac{\eps}{2},
\end{equation*}
so any cover of $(A \cap W_i) \setminus R$ by $\delta$-balls contains $\gtrsim \delta^{\eps - \Delta}$ balls.

\vs{0.3}

Now, \eqref{eq:plane-Hausdorff-content} and \eqref{eq:square-Hausdorff-content} also entail that there exist distinct $R_1, ..., R_{2^{n-k} + 1}$ whose intersections with $A \cap W_i$ each have positive $(\Delta - \eps)$-dimensional Hausdorff content. In particular,
\begin{equation} \label{eq:checker-square-covers}
    \card \big\{ x \in A' \!: B(x,\delta) \cap (A \cap W_i \cap R_p) \neq \varnothing \big\} \gtrsim \delta^{\eps - \Delta}
\end{equation}
for $p = 1, ..., 2^{n-k} + 1$ because $\{ B(x,\delta) \!: x \in A' \}$ is a cover of $A$. Necessarily, at least $2$ of these squares $R_p, R_q$ are mutually non-adjacent, so they are $4\eta$-separated. Therefore, if $x,y \in A'$ are such that
\begin{gather*}
    B(x,\delta) \cap (A \cap W_i \cap R_p) \neq \varnothing \quad \text{and} \quad B(y,\delta) \cap (A \cap W_i \cap R_q) \neq \varnothing, \\
    \intertext{then}
    \| x - y \| > 4\eta - 2\delta > 2\eta,
\end{gather*}
so that $x \backsim_{\scriptscriptstyle{V}} y$. In conjunction with \eqref{eq:checker-square-covers}, this yields the estimate
\begin{align*}
    \card \big\{ & (x,y) \in (A')^2 \!:  B(x,\delta) \cap (A \cap W_i \cap R_p) \neq \varnothing, \\
    &\hs{1.1} B(y,\delta) \cap (A \cap W_i \cap R_q) \neq \varnothing, \, x \backsim_{\scriptscriptstyle{V}} y \big\} \gtrsim \delta^{2(\eps - \Delta)}.
\end{align*}
No $\delta$-ball intersects more than $3^k$ fat planes in $\mathcal{T}_V$, so we may sum the previous over all $i \in \{ 1, ..., j \}$ to get
\begin{equation*}
    \card \big\{ (x,y) \in (A')^2 \!: x \backsim_{\scriptscriptstyle{V}} y \big\} \gtrsim j (\delta^{\eps - \Delta})^2 \gtrsim \delta^{(\Delta - \gamma + \eps) + 2(\eps - \Delta)} \geq \delta^{3\eps + s - 2\gamma},
\end{equation*}
where the final inequality follows from our original hypothesis that $\Delta \geq \gamma - s$.  This is the desired lower bound on the individual summands in \eqref{eq:energy}, and multiplying by $\card E'$ yields the desired bound on $\mathcal{E}$ itself:
\begin{equation*}
    \mathcal{E} \gtrsim \delta^{3\eps + s - 2\gamma} \card E'.
\end{equation*}
\textsc{Step 4.} In combination with our upper bound $\delta^{-k(n-k)+k-2\gamma'} \gtrsim \mathcal{E}$, this at last provides a concrete upper bound on $\card E'$ in terms of $\delta$, namely,
\begin{equation*}
    \card E' \lesssim \delta^{-k(n-k)+k - 2\gamma' - (3\eps + s - 2\gamma)} = \delta^{-k(n-k) + (k-s) - 3\eps - 2(\gamma' - \gamma)}.
\end{equation*}
Since $\gamma' > \gamma$ was arbitrary and does not appear in the implicit constant, the estimate $\card E' \lesssim \delta^{-k(n-k) + (k-s) - 3\eps}$ follows at once. This holds for every $\delta$-separated subset $E' \subseteq E_{s,\eps}$, so we conclude \eqref{eq:main-reduction} and, in turn, \eqref{eq:second-DC-lem}. \hfill $\square$

\hypertarget{appendix}{}
\section*{Appendix}

This section details the proof of Proposition \ref{prop:rams} from Theorem \ref{thm:rams}. We make only minor modifications to the notation of \cite{rams2002packing}, to which the reader is referred for thorough definitions of the pertinent concepts. In summary:
\begin{enumerate}[label=\textbullet, noitemsep, topsep=-3pt]
    \item $\rho_e \!: \R^n \to \R^{n-1}$ is the orthogonal projection onto $(\Span e)^\perp \cong \R^{n-1}$.
    \item Given $N \in \Z_+$, $\Sigma := \{ 1, ..., N \}^\N$ denotes our symbol space, and a typical element of $\Sigma$ is written as $\omega = (\omega_1, \omega_2, ...)$. We also denote $\omega^m := (\omega_1, ..., \omega_m)$.
    \item Given an IFS $(g_i)_{i=1}^N$, $\Pi \!: \Sigma \to \R^{n-1}$ indicates the projection map
    \begin{equation*}
        \Pi(\omega) := \lim_{m \to \infty} g_{\omega^m}(0),
    \end{equation*}
    where $g_{\omega^m} := g_{\omega_1} \circ \cdots \circ g_{\omega_m}$. This limit always exists by the Cantor intersection theorem. When working with a parametrized family $\big( f_i(\,\cdot\,;t) \big)_{i=1}^N$, the projections $\Pi_t$ are distinguished with a subscript $t$.
\end{enumerate}

Lastly, we restate Definition 4.4 of \cite{rams2002packing}. A smoothly parametrized family $\big( f_i(\,\cdot\,;t) \big)_{i=1}^N$ as above is said to be \textit{transverse} (or to satisfy the \textit{transversality condition}) if there exists a constant $L > 0$ with the following property: for each parameter value $u$ and each pair $\omega, \kappa \in \Sigma$ with $\omega_1 \neq \kappa_1$,
\begin{equation} \label{eq:transversality}
    \| \Pi_u(\omega) - \Pi_u(\kappa) \| < L \quad \text{implies} \quad \Big| \det D_t \big( \Pi_t(\omega) - \Pi_t(\kappa) \big) \big|_{t=u} \Big| > L.
\end{equation}
Here, $D_t$ denotes the total derivative with respect to $t$. Loosely, this definition states that, whenever the points $\Pi_t(\omega), \Pi_t(\kappa)$ are close for some $t=u$, they do not remain close for long as the parameter $t$ changes. The condition $\omega_1 \neq \kappa_1$ suggests that $\Pi_t(\omega)$ and $\Pi_t(\kappa)$ ``should" land in different regions of the attractor $K_t$ for any value of $t$, although these regions may overlap for values of $t$ at which the IFS is degenerate.

\define{Proof of Proposition \ref{prop:rams}.} It suffices to work in local coordinates, so we let $V \subset \bbs^{n-1}$ be an open set whose closure $\overline{V}$ is diffeomorphic to a bounded subset of $\R^{n-1}$. These local coordinates also afford us a consistent identification of the tangent hyperplanes to $\bbs^{n-1}$ with $\R^{n-1}$. Dispensing with these technicalities, we we simply refer to our parameter space as $\bbs^{n-1}$ and use the formula $\rho_e(x) = x - (x \cdot e) e$ for the orthogonal projections.

\textsc{Step 1.} Let $(g_i)_{i=1}^N$ be an IFS on $\R^n$ with limit set $K$ and satisfying the SSC. We seek to produce a smooth family of IFS on $\R^{n-1}$ to which we can apply Theorem 1.1 of \cite{rams2002packing}.

To this end, we define $\big( f_i(\,\cdot\,;e) \big)_{i=1}^N$ by
\begin{equation*}
    f_i(\xi;e) := (\rho_e \circ g_i)\big( \rho_e^{-1}(\xi) \big)
\end{equation*}
for each $e \in \bbs^{n-1}$, where $\rho_e^{-1}(\xi)$ is any preimage of the point $\xi \in \R^{n-1}$. This definition is unambiguous because $g_i$ takes the form $g_i(x) = ax + b$ for some $a \in \R$ and $b \in \R^n$, whence
\begin{equation*}
    (\rho_e \circ g_i)\big( \rho_e^{-1}(\xi) \big) = \rho_e\big( a \rho_e^{-1}(\xi) + b \big) = a \rho_e\big( \rho_e^{-1}(\xi) \big) + \rho_e(b) = a\xi + \rho_e(b)
\end{equation*}
for any choice of $\rho_e^{-1}(\xi)$. This also shows that $\big( f_i(\,\cdot\,;e) \big)_{i=1}^N$ is smooth in both $x$ and $e$, as $\rho_e(b) = b - (b \cdot e) e$.

\textsc{Step 2.} We show that $\big( f_i(\,\cdot\,;e) \big)_{i=1}^N$ is a transverse family. Let $\omega, \kappa \in \Sigma$, where $\omega_1 \neq \kappa_1$, and let $f_{\omega^m}(\xi;e)$ denote the composite map
\begin{equation*}
    f_{\omega_1}(\,\cdot\,;e) \circ \cdots \circ f_{\omega_m}(\,\cdot\,;e)
\end{equation*}
evaluated at $\xi$. Then
\begin{align*}
    f_{\omega^m}(\xi;e) &= \big( (\rho_e \circ g_{\omega_1} \circ \rho_e^{-1}) \circ \cdots \circ (\rho_e \circ g_{\omega_m} \circ \rho_e^{-1}) \big)(\xi) \\
    &= \big( \rho_e \circ (g_{\omega_1} \circ \cdots \circ g_{\omega_m}) \circ \rho_e^{-1} \big)(\xi) \\
    &= (\rho_e \circ g_{\omega^m} \circ \rho_e^{-1})(\xi)
\end{align*}
for any section $\rho_e^{-1}$ of $\rho_e$. Therefore, by the continuity of $\rho_e$,
\begin{equation*}
    f_\omega(\xi;e) := \lim_{m \to \infty} f_{\omega^m}(\xi;e) = \rho_e\!\left( \lim_{m \to \infty} g_{\omega^m} \big( \rho_e^{-1}(\xi) \big) \right) = (\rho_e \circ g_\omega)\big( \rho_e^{-1}(\xi) \big).
\end{equation*}
In particular, we can take $\rho_e^{-1}(0) = 0$, so
\begin{equation*}
    \Pi_e(\omega) = f_\omega(0;e) = (\rho_e \circ g_\omega)(0) = \rho_e(\Pi(\omega));
\end{equation*}
likewise for $\kappa$.

Denote $z = (z_1, ..., z_n) = \Pi(\omega) - \Pi(\kappa)$, and suppose
\begin{equation} \label{eq:transversality-hyp}
    \| \rho_u(z) \| = \big\| \rho_u(\Pi(\omega)) - \rho_u(\Pi(\kappa)) \big\| = \| \Pi_u(\omega) - \Pi_u(\kappa) \| < \frac{c}{\sqrt{2}}
\end{equation}
for some $u \in \bbs^{n-1}$, where $c \in (0,1]$ is a constant such that $\dist(g_i(K),g_j(K)) > c$ for all $i \neq j$. Such an $c$ exists because $(g_i)_{i=1}^N$ satisfies the SSC. Since $\Pi(\omega) \in g_{\omega_1}(K)$, $\Pi(\kappa) \in g_{\kappa_1}(K)$, and $\omega_1 \neq \kappa_1$, it follows that $\|z\| > c$\textemdash a fact we shall use shortly.

To show transversality, we must compute
\begin{equation*}
    \det D_e\big( \rho_e(z) \big) \big|_{e=u}.
\end{equation*}
The determinant is invariant under a linear change of coordinates, so we can rotate our coordinate system so that $u = e_n = (0, ..., 0, 1)$. Consider $h \!: e \mapsto \rho_e(z)$ as a map from $\R^n$ to $\R^n$, i.e., by extending $\rho_e(z)$ to take parameter values in $\R^n$. Considered as an $n \times n$ matrix, the $j$th column of the derivative $D_e h(e) |_{e=e_n}$ is given by the directional derivative
\begin{align*}
    \frac{d}{dr} \rho_{e_n + re_j}(z) \big|_{r=0} &= \frac{d}{dr} \big( z - (z \cdot (e_n + re_j)) \+ (e_n + re_j) \big) \big|_{r=0} \\
    &= -z_j (e_n + re_j) - (z_n + rz_j) e_j \big|_{r=0} = -z_j e_n - z_n e_j,
\end{align*}
yielding
\begin{equation*}
    D_e h(e) |_{e=e_n} = \left( \begin{array}{cccc}
    -z_n & \cdots & 0 & 0 \\
    \vdots & \ddots & \vdots & \vdots \\
    0 & \cdots & -z_n & 0 \\
    -z_1 & -z_2 & \cdots & -2z_n
    \end{array} \right).
\end{equation*}
Since $D_e h(e) |_{e=e_n}$ restricts to an automorphism of the tangent plane $(\Span e_n)^\perp \cong T_{e_n} \bbs^{n-1}$, and since the standard coordinate frame $(e_1, ..., e_n)$ is adapted to $\bbs^{n-1}$ at the north pole $e_n$, the matrix of this restricted linear map is obtained simply by omitting the $n$th row and $n$th column of the matrix. That is, $D_e\big( \rho_e(z) \big) \big|_{e=e_n} = -z_n I_{n-1}$ and, consequently,
\begin{equation*}
    \det D_e\big( \rho_e(z) \big) \big|_{e=e_n} = \det\!\big( \!-\!z_n I_{n-1} \big) = (-z_n)^{n-1}.
\end{equation*}
Now, since $z = \rho_{e_n}(z) + z_n e_n$, $\|\rho_{e_n}(z)\|^2 < 2^{-1} c^2$, and $\|z\|^2 > c^2$, it must be that $|z_n|^2 = \| z_n e_n \|^2 > 2^{-1} c^2$ and, in turn,
\begin{equation*}
    \Big\| D_e\big( \rho_e(z) \big) \big|_{e=e_n} \Big\| = \| z_n \|^{n-1} > \frac{c^{n-1}}{2^{(n-1)/2}} \geq \frac{c}{\sqrt{2}}.
\end{equation*}
In view of Equation \eqref{eq:transversality-hyp}, we conclude that, whenever $\omega_1 \neq \kappa_1$,
\begin{equation*}
    \big\| \Pi_u(\omega) - \Pi_u(\kappa) \big\| < \frac{c}{\sqrt{2}} \quad \text{implies} \quad \Big| \det D_e\big( \Pi_e(\omega) - \Pi_e(\kappa) \big) \big|_{e=u} \Big| = \Big\| D_e\big( \rho_e(z) \big) \big|_{e=u} \Big\| > \frac{c}{\sqrt{2}},
\end{equation*}
so $\big( f_i(\,\cdot\,;e) \big)_{i=1}^N$ satisfies the transversality condition.

\textsc{Step 3.} We apply Theorem \ref{thm:rams} to get
\begin{equation*}
    \dim_P \, \{ e \in \bbs^{n-1} \!: \dim K_e \leq s \} \leq s \qquad \forall \, 0 \leq s < \min \left\{ \, n-1, \, \sup_{e \in \bbs^{n-1}} \+ \sigma(e) \right\}.
\end{equation*}
If $\dim K \geq n-1$, then $\sup_{e \in \bbs^{n-1}} \+ \sigma(e) = n-1$ and the desired conclusion holds for all $0 \leq s < \dim K$, the values of $s$ in the interval $(n-1,\dim K)$ giving a trivial bound. If instead $\dim K < n-1$, then $\sup_{e \in \bbs^{n-1}} \+ \sigma(e) = \dim K$, and again the bound on the exceptional set holds for all $0 \leq s < \dim K$. \hfill $\square$

\phantomsection
\section*{Acknowledgement}
My sincere thanks go to my advisor, Bobby Wilson, for his guidance, encouragement, and feedback throughout the writing of this paper. I also extend my gratitude to the referees for their kind and helpful suggestions, including an improved statement of the main theorem.

\bibliographystyle{plain}
\bibliography{references}

\end{document}